\documentclass[conference]{IEEEtran}
\IEEEoverridecommandlockouts
\usepackage{tcolorbox}
\newcommand{\myexample}[2]{
	\begin{tcolorbox}[colback=black!5!white,colframe=black,title={Identity #1}]
		#2
	\end{tcolorbox}
}

\newcommand{\myexamplee}[2]{
	\begin{tcolorbox}[colback=black!5!white,colframe=black,title={  #1}]
		#2
	\end{tcolorbox}
}

\usepackage{cite}
\usepackage{amsmath,amssymb,amsfonts}
\usepackage{algorithmic}
\usepackage{graphicx}
\usepackage{textcomp}
\usepackage{xcolor} 
\usepackage{epsfig}
\usepackage[font=small,skip=0pt]{caption}
\usepackage{amsmath,amssymb,amsthm}
\usepackage{comment} 
\usepackage{tabularx}
\usepackage{algorithmic}
\usepackage{graphicx}
\usepackage{algorithm}
\usepackage{enumitem}
\usepackage{eufrak}
\usepackage{lipsum}
\usepackage{stfloats}
\usepackage{bbm}
\usepackage{relsize}
\usepackage{multicol}

\newtheorem{thm}{Theorem}
\newtheorem{lem}{Lemma}
\newtheorem{remark}{Remark}

\def\bsum{\mathlarger{\sum}}

\def\BibTeX{{\rm B\kern-.05em{\sc i\kern-.025em b}\kern-.08em
		T\kern-.1667em\lower.7ex\hbox{E}\kern-.125emX}}
\begin{document}
	
	\title{Two identities for Poisson Point Processes and Voronoi Tessellations with Applications}
	
	\author{
		\IEEEauthorblockN{Mohsen Amidzadeh}\vspace{-0 pt}\\
		\IEEEauthorblockA{Department of Computer Science, Aalto University, Finland\vspace{-0 pt}}
}
\maketitle

%----------------------Abstract------------------
\begin{abstract}
In this paper, we introduce two identities—one pertaining to the state space of Poisson Point Processes (PPPs), 
and the other for the Voronoi tessellations formed by PPPs. 
Then, we explore several applications of these identities 
within the context of wireless cellular networks.
\end{abstract}

%\newpage
\begin{IEEEkeywords}
	Poisson Point Process, Vornoi Tessellation, Cellular Networks.
\end{IEEEkeywords}

%/\/\/\/\/\/\/\/\/\/\/\/\/\/\/\/\/\/\/\/\/\/\/\/\/\/\/\/\/\/\/\/\/\/\/\/\/\/\/\/\/\/\/\/\/\/\/\/\/\/\/\/\/\/\/\/%
%  													   	 Introductoin    										%
%/\/\/\/\/\/\/\/\/\/\/\/\/\/\/\/\/\/\/\/\/\/\/\/\/\/\/\/\/\/\/\/\/\/\/\/\/\/\/\/\/\/\/\/\/\/\/\/\/\/\/\/\/\/\/\/%
\section{Introduction}\label{Sec_Intro}
Stochastic geometry is considered a versatile tool for modeling and analyzing wireless cellular networks, 
both for homogeneous and Heterogeneous Networks (HetNets) \cite{Andrews2011,Baccelli2009,Baccelli2009_vol2,Dhillon2012,Haenggi2009,Han_Shin2012,Andrews2016,Ours_2020,Ours_2021,Ours_2022_twc,Ours_2022}. 
In the context of caching at the wireless edge, Poisson Point Processes (PPPs) have been used to model deployment of BSs and locations of UE in \cite{Lee2015,Wu2017,Ye2019,Xu2017,Choi2021}. 
In \cite{Lee2015}, where BSs apply a beamforming content delivery and the expected ergodic spectral efficiency is optimized, stochastic geometry is used to model BS deployment.
In  \cite{Serbetci2017}, optimum cache policies are found in HetNets 
modelled based on independent PPPs.
%to achieve an optimum cache policy based on total miss probability perspective.
In \cite{Wu2017}, two independent PPPs are exploited to model the deployment of two tiers of a HetNets. 
Aiming to find an optimal cache placement policy, independent PPPs are exploited in
\cite{Ye2019} to obtain an expression for the successful transmission probability in a HetNet.
In \cite{Xu2017}, the coverage probability is approximated using PPP for the multi-antenna small-cell networks to design an optimal cache placement policy. 
In \cite{Choi2021}, the authors leverage PPP to optimize the minimum of the cache hit rates of
different request-related categorized UEs.
\cite{Ours_2020,Ours_2021,Ours_2022_twc} apply stochastic geometry to demonstrate that the probabilistic cache policies significantly improves the quality-of-service.
\cite{Ours_2022} presents a compound cache policy 
schemes with performance analysis done by the aid of stochastic geometry. 
These studies highlight the applicability of stochastic geometry in analysis the performance metrics 
of cellular networks.
In this paper we present two identities for the Voronoi tessellations constructed of Poisson Point Processes,
and provide some applications of these identities in the context of mobile cellular networks.
This work provide a versatile analytical tool for analyzing the conventional transmission schemes
of the wireless networks using the stochastic geometry.

The remainder of this paper is organized as follows. 
In Section \ref{Sec_Identity}, we present two identities related to Poisson point processes.
In Section~\ref{Sec_Application}, we discuss about the applications of these identities.
Section \ref{Sec_Conclu} then concludes the paper.

\textbf{Notations}: We use lower-case $a$ for scalars,
bold-face lower-case $\textbf{a}$ for vectors and bold-face uppercase
$\textbf{A}$ for matrices. 
%$\textbf{A}^\top$ and $\textbf{A}^\dagger$ are the transpose and Hermitian conjugate of
%$\textbf{A}$ and 
$\|\textbf{a}\|$ is the Euclidean norm of $\textbf{a}$,
$\{a_n\}_n$ collects the components of vector $\textbf{a}$
and $\mathbbm{1}(\cdot)$ is the indicator function.
%and $\{a_n\}_{n=1}^N$ collects the components of vector $\boldsymbol{a}$
%from $n=1$ to $n=N$. 
%The notation $[ a_{i} ]_{i} $ shows a column vector with i-th component being $a_i$.
%The notation  $\textbf{a}^{-1}$ denotes the component-wise inversion of  $\textbf{a}$,
%and $\textbf{1}$ and $\textbf{0}$ to denote
%the vector with all elements equal to one and zero, respectively.
%Further, $\delta_{ij}$ is the Kronecker delta function, 
%$\mathbb{C}$ and $\mathbb{R}$ denote the complex valued and real valued numbers.

%/\/\/\/\/\/\/\/\/\/\/\/\/\/\/\/\/\/\/\/\/\/\/\/\/\/\/\/\/\/\/\/\/\/\/\/\/\/\/\/\/\/\/\/\/\/\/\/\/\/\/\/\/\/\/\/%
%  													   	 System Model  											%
%/\/\/\/\/\/\/\/\/\/\/\/\/\/\/\/\/\/\/\/\/\/\/\/\/\/\/\/\/\/\/\/\/\/\/\/\/\/\/\/\/\/\/\/\/\/\/\/\/\/\/\/\/\/\/\/%
\section{An identity for Poisson-Voronoi Tessellations} \label{Sec_Identity}
We now present the first identity in the following Lemma.
\begin{lem}\label{Lemma1}
	Let $\{\boldsymbol{x}_i\}_{i\in\Phi_1}$ be the points of a homogeneous PPP $\Phi_1$ with intensity $\lambda_1$,
	and $S(\cdot)$ and $P(\cdot)$ be two real-valued functions on the state space of $\Phi_1$. We then have:
	\begingroup\makeatletter\def\f@size{9}\check@mathfonts
	\myexample{1: Expectation of sum-product of functions on the state space of PPP.}{\vspace{-10pt}	
	\begin{align*}
		&\mathbb{E}\left\{ \prod_{k\in \Phi_1} P(\boldsymbol{x}_k) \sum_{k\in \Phi_1} S(\boldsymbol{x}_k) \right \} 
		= \\ &		~~~~~~~~~
		\lambda_1 \iint_{\mathbb{R}^2} S(\boldsymbol{s})P(\boldsymbol{s}) d\boldsymbol{s} \: \exp\left( \lambda_1 \iint_{\mathbb{R}^2} \big( P(\boldsymbol{s}) - 1\big) d\boldsymbol{s} \right)
	\end{align*}\vspace{-10pt}	
	}
	\endgroup     
\end{lem}
\begin{proof}
please refer to Appendix~\ref{App1} for the proof.
\end{proof}
The second identity is presented as follows.
\begin{lem}\label{Lemma2}
		Assume $\{\boldsymbol{r}_i\}_{i\in \Phi_2}$ to be the Poisson points located based on a PPP $\Phi_2$,
		and  consider Voronoi tessellation constructed by these points 
		with $\mathcal{V}_0$ being the Voronoi cell associated with $\boldsymbol{r}_0$,
		then the quantity $\mathbbm{1}(\boldsymbol{x} \in \mathcal{V}_0)$, for any location $\boldsymbol{x}\in\mathbb{R}^2$,
		can be expressed based on the following function product in $\Phi_2$:
	\begingroup\makeatletter\def\f@size{9.5}\check@mathfonts
	\myexample{2: Expression of $\mathbbm{1}(\boldsymbol{x}\in \mathcal{V}_0)$ based on a  product of functions.}{\vspace{-10pt}	
	\begin{align}\label{EQ:State1}
		\mathbbm{1}(\boldsymbol{x} \in \mathcal{V}_0) 
		= \prod_{k\in \Phi_2 \backslash \{0\} } \mathbbm{1}\big( \boldsymbol{x} \in M_0\left(\boldsymbol{r}_k\right) \big),
	\end{align}%\vspace{-10pt}	
	}
	\endgroup
	\noindent where $\mathbbm{1}(\cdot)$ is the indicator function and 
	\begin{align}\label{EQ:State2}
			M_0(\boldsymbol{r}_k) := \big\{\boldsymbol{r}'\in \mathbb{R}^2 ~\big|~ \|\boldsymbol{r}'-\boldsymbol{r}_0\| \leq \|\boldsymbol{r}'-\boldsymbol{r}_k\| \big\}.
	\end{align} 
	Moreover, for the space $\big\{ \boldsymbol{r} \in \mathbb{R}^2 \: |\: \textbf{x} \in M_0(\boldsymbol{r}) \big\}$ with a given $\textbf{x} \in \mathbb{R}^2$, 
	we have:
	\begingroup\makeatletter\def\f@size{9.5}\check@mathfonts
	\myexamplee{}{\vspace{-4pt}	
	$$
	\big\{ \boldsymbol{r} \in \mathbb{R}^2 \: |\: \textbf{x} \in M_0(\boldsymbol{r}) \big\} = B_c(\boldsymbol{r};\boldsymbol{x},\|\boldsymbol{x}\|),
	$$	
	}
	\endgroup
	\noindent where 
	$B_c(\boldsymbol{r};\boldsymbol{x},\|\boldsymbol{x}\|)$ defines the space exterior to the circle centered at $\boldsymbol{x}$ 
	with a radius of $\|\boldsymbol{x}\|$.
%	\begingroup\makeatletter\def\f@size{8.}\check@mathfonts
%	\begin{align*} 
%	&\mathbb{E}\left\{ \bsum_{i \in \Phi_1,\:\textbf{x}_i\in\mathcal{V}_0} \bsum_{k \in \Phi_2\backslash\{0\}}\!\!\!f(\textbf{x}_i,\textbf{r}_k) \right \} 
%	= 	\mathbb{E}\left\{ \bsum_{i \in \Phi_1}  \bsum_{k \in \Phi_2\backslash\{0\}}\!\!\!f(\textbf{x}_i,\textbf{r}_k) \mathbbm{1}(\textbf{x}_i\in\mathcal{V}_0)\right \} \\
%	&= \mathbb{E}_{\Phi_1,\Phi_2}\left\{ \bsum_{i \in \Phi_1}  \bsum_{k \in \Phi_2\backslash\{0\}}\!\!\!f(\textbf{x}_i,\textbf{r}_k) \prod_{k\in \Phi_2 \backslash \{0\} } \mathbbm{1}\big( \boldsymbol{x}_i \in M_0\left(\boldsymbol{r}_k\right) \big) \right\} \\
%	&= \mathbb{E}_{\Phi_1}\Bigg\{ \sum_{i\in\Phi_1}\lambda_2\iint_{\mathbb{R}^2} f(\textbf{x}_i,\textbf{r}) \mathbbm{1}(\textbf{x}_i\in M_0(\textbf{r}))rdrd\theta\\
%	&\qquad\qquad\times
%	\exp\left( \lambda_2 \iint_{\mathbb{R}^2}\big( \mathbbm{1}(\textbf{x}_i\in M_0(\textbf{r}))-1 \big)rdrd\theta \right)\Bigg\} \\
%	&= \mathbb{E}_{\Phi_1}\Bigg\{ \sum_{i\in\Phi_1}\lambda_2\iint_{B_c(\textbf{r};\textbf{x}_i,\|\textbf{x}_i\|)} \hspace{-25 pt}f(\textbf{x}_i,\textbf{r}) rdrd\theta\:
%	\exp\left( -\lambda_2 \iint_{B_c^\prime(\textbf{r};\textbf{x}_i,\|\textbf{x}_i\|)}\hspace{-25 pt}1\:rdrd\theta \right) \Bigg\} \\
%	&= \lambda_1 \iint_{\mathbb{R}^2}\lambda_2\iint_{B_c(\textbf{r};\textbf{x},\|\textbf{x}\|)} \hspace{-30pt}f(\textbf{x},\textbf{r}) rdrd\theta\:
%	\exp\left( -\lambda_2 \pi x^2\right)  xdxd\theta \\
%	\end{align*} 
%	\endgroup
%	for the function $f(\cdot):\mathbb{R}^2 \to \mathbb{R}$.
\end{lem}
\begin{proof}
please refer to Appendix~\ref{App2} for the proof.
\end{proof}

\begin{remark}\label{Remark}
	Consider $\{\textbf{x}_i\}_{i\in\Phi_1}$ and $\{\textbf{r}_i\}_{i\in\Phi_2}$
	as Poisson points of two independent PPPs $\Phi_1$ and $\phi_2$ with intensities $\lambda_1$ and $\lambda_2$,
	respectively, and $\mathcal{V}_0$ being the Voronoi cell associated with $\boldsymbol{r}_0$,
	and $P(\cdot)$ a real-valued function on the state space of $\Phi_1$.
	Then, we have:
	\begingroup\makeatletter\def\f@size{8.5}\check@mathfonts
	\begin{align*} 
	&\mathbb{E}\left\{ \sum_{k\in \mathcal{V}_0} P(\textbf{x}_k) \right\} 
	= \lambda_1 \iint_{\mathbb{R}^2} P(\textbf{x})\exp\Big(-\pi\lambda_2  \|\textbf{x}\|^2 \Big) d\textbf{x},
	\end{align*}
	\endgroup
	and
	\begingroup\makeatletter\def\f@size{8.5}\check@mathfonts
	\begin{align*} 
	&\mathbb{E}\left\{ \prod_{k\in \mathcal{V}_0} P(\textbf{x}_k) \right\} 
	= \exp\left( \lambda_1 \!\!\iint_{\mathbb{R}^2} \!\!\Big(P(\textbf{x})\exp\Big(-\pi\lambda_2  \|\textbf{x}\|^2 \Big) -1 \Big) d\textbf{x} \right).
	\end{align*}
	\endgroup	
\end{remark}
\begin{proof}
	See Appendix \ref{App3} for the proof.
\end{proof}
%/\/\/\/\/\/\/\/\/\/\/\/\/\/\/\/\/\/\/\/\/\/\/\/\/\/\/\/\/\/\/\/\/\/\/\/\/\/\/\/\/\/\/\/\/\/\/\/\/\/\/\/\/\/\/\/%
%  													   	 SubSec_Placement										%
%/\/\/\/\/\/\/\/\/\/\/\/\/\/\/\/\/\/\/\/\/\/\/\/\/\/\/\/\/\/\/\/\/\/\/\/\/\/\/\/\/\/\/\/\/\/\/\/\/\/\/\/\/\/\/\/%
\section{Applications} \label{Sec_Application}
In the following sections, we present two applications of developed identities
for analyzing network performance within cellular networks.

\subsection{Average Number of UEs within a Voronoi Cell}\label{Subsec:app1}

\textcolor{black}{We here intend to use Lemma~\ref{Lemma2}
to compute the average number of UEs within a Voronoi Cell of BSs.
This is an important parameter since 
each UE is served by its nearest BS
in the conventional transmission scheme of cellular network.}

For this, consider BSs and UEs which are located  based on two independent PPPs, 
$\Phi_u$ and $\Phi_b$ with intensities $\lambda_u$ and $\lambda_b$.
We denote the location of BSs and UEs by $\{\textbf{r}_i\}_{i \in \Phi_b}$ and $\{\textbf{x}_i\}_{i\in\Phi_u}$, 
respectively.
The BSs constitutes a Voronoi tessellation with different cells denoted by $\{\mathcal{V}_i\}_i$.
Notice that the PDF of cell size of Voronoi tessellation
is accurately \textit{approximated} by a gamma distribution with shape $K$ and scale $1/(K\lambda_b)$ \cite{Pineda2007,Cao2012} with $K=3.575$.
This gives the average cell size equal to $\bar{A} = K \frac{1}{K\lambda_b} = \frac{1}{\lambda_b}$.
Therefore, the average number of UEs within a typical Voronoi cell $\bar{n}_u$ 
can be \textit{approximated} as follows:
\begin{align}\label{EQ:AverageN}
\bar{n}_u = \lambda_u \: \bar{A} = \frac{\lambda_u }{\lambda_b}.
\end{align}
We now intend to verify this quantity using Lemma~\ref{Lemma2}.
Without loss of generality, we can compute $\bar{n}_u$ for $\mathcal{V}_0$.
\begingroup\makeatletter\def\f@size{8.75}\check@mathfonts
\begin{align}\label{EQ:avg_n}
\bar{n}_u &= \mathbb{E}\left\{ \sum_{i \in \Phi_u} \mathbbm{1}(\textbf{x}_i \in \mathcal{V}_0) \right\} \notag  \\
&\overset{(a)}= \mathbb{E}_{\Phi_b}\mathbb{E}_{\Phi_u}\left\{ \sum_{i \in \Phi_u} \prod_{k \in \Phi_b\backslash\{0\}} \mathbbm{1}(\textbf{x}_i \in M_0(\textbf{r}_k))  \right\}  \notag \\
&\overset{(b)}= \mathbb{E}_{\Phi_u}\left\{ \sum_{i\in\Phi_u} 
\exp\left( \lambda_b \iint_{\mathbb{R}^2}\big( \mathbbm{1}(\textbf{x}_i\in M_0(\textbf{r}))-1 \big)d\textbf{r} \right) \big| \Phi_b\right\}  \notag \\
&\overset{(c)}= \mathbb{E}_{\Phi_u}\left\{ \sum_{i\in \Phi_u} \exp\left( \lambda_b \iint_{B_c(\boldsymbol{r};\textbf{x}_i,\|\textbf{x}_i\|)} \hspace{-31pt}1\: d\textbf{r} - \lambda_b\iint_{\mathbb{R}^2}\hspace{-5pt}1\: d\textbf{r} \right) \big| \Phi_b \right\}  \notag \\
&\overset{(d)}= \lambda_u \iint_{\mathbb{R}^2} \exp\left( \lambda_b \iint_{B_c(\boldsymbol{r};\textbf{x},\|\textbf{x}\|)} \hspace{-25pt}1\: d\textbf{r} - \lambda_b\iint_{\mathbb{R}^2}\hspace{-5pt}1\: d\textbf{r} \right) xdx\:d\theta  \notag \\
&= \lambda_u \iint_{\mathbb{R}^2} \exp\left( -\lambda_b \pi x^2 \right) xdx\:d\theta \:=\: \frac{\lambda_u}{\lambda_b},
\end{align}
\endgroup
where $x=\|\textbf{x}\|$ and $r=\|\textbf{r}\|$.
Furthermore, for (a) we used the Lemma~\ref{Lemma2}, for (b) we leveraged Campbell’s Theorem \cite{Kingman93}  for $\Phi_b$,
and for (c), we used 
$$
\big\{ \boldsymbol{r} \in \mathbb{R}^2 \: |\: \textbf{x} \in M_0(\boldsymbol{r}) \big\} = B_c(\boldsymbol{r};\boldsymbol{x},\|\boldsymbol{x}\|),
$$
from Lemma~\ref{Lemma2} with $B_c(\boldsymbol{r};\boldsymbol{x},\|\boldsymbol{x}\|)$ 
denoting the space outside 
the circle centered at $\boldsymbol{x}$ with a radius of $\|\boldsymbol{x}\|$.
Moreover, for (d), Campbell’s Theorem was now used for $\Phi_u$.
The obtained result interestingly confirms the approximation obtained in Eq. \eqref{EQ:AverageN}.

\subsection{Total Bandwidth Consumption}\label{Subsec:app2}
In this section we use Lemmas~\ref{Lemma1} and \ref{Lemma2} to compute the total bandwidth consumption
of the conventional unicast scheme \cite{Wen2017, Li2018, Wu2020} within the cellular networks.
For this, consider a network with BSs applying the conventional single-point unicast scheme 
to serve the requesting UEs of the network.
Each UE of this transmission scheme requests the needed content from its nearest BS.
We consider the transmission rate $R$ for the content requests of UEs.
Therefore, this scheme develops an on-demand transmission with BS-specific association.
Each BS uses the full bandwidth $W$ for content delivery, i.e., frequency reuse 1 is used \cite{rusek2013}.
This is according to the concurrent standards of mobile networks~\cite{laiho2006,garcia2009,rusek2013}.

The locations of UEs and BSs are considered based on two independent PPPs, 
$\Phi_u$ and $\Phi_b$, respectively, and with intensities 
$\lambda_u$ and $\lambda_b$.

Since each UE is served by its nearest BS, they establish a Voronoi tessellation with different cell size.
As such, there may be  $U$ distinct UEs requesting from the BS of a typical Voronoi cell.
For UE $k$ served by the network, the SINR is:
\begingroup\makeatletter\def\f@size{10.}\check@mathfonts
\begin{align}\label{EQ_SINR}
\gamma_{k}=\frac{ |h_0^k|^2 \|\boldsymbol{x}_k - \boldsymbol{r}_0\|^{-e}}{ {1}/{\gamma_{\rm tx}} 
+ \underbrace{\sum_{j\in\Phi_b \backslash \{0\}} |h_j^k|^2 \|\boldsymbol{x}_k - \boldsymbol{r}_j\|^{-e}}_{I_k}},
\end{align}
\endgroup
where the  Tx-SNR is $\gamma_{\rm tx}^{\rm UC}=p_{\rm tx}/(WN_0)$
with $p_{\rm tx}$ being the average transmission power of BSs.
Moreover, $e$ is the path-loss component, 
$h_0^k$ is the channel coefficient between nearest BS and UE $k$, 
and $h_j^k $ is the channel coefficient from BS $j$ to UE $k$. 
We use a standard distance-dependent path-loss model, and assume a
Rayleigh distribution for the channel coefficient, i.e.,
$|h_{j,k}|^2\sim \exp(1)$.

Within each cell, the BS follows a UE-specific resource allocation policy
to ensure an adequate allotment of resources for the served UEs. 
However,
when the SINR of a UE diminishes significantly, near infinite bandwidth is
required to serve the UE. 
To address this challenge, it is conventional \cite{goldsmith2005} to implement a truncated SINR strategy,
wherein the bandwidth allocated to serve UE $k$ is:
\begingroup\makeatletter\def\f@size{10.}\check@mathfonts
\begin{align}\label{EQ_serviceTh}
w_k^{\rm UC} = \begin{cases}
	{R}/{\log_2(1+\gamma_{k}^{\rm UC})}, & \gamma_{k}^{\rm UC} \geq \gamma\\
	0,  & {\rm otherwise}
\end{cases},
\end{align}
\endgroup
where $\gamma$ is the SINR threshold and $R$ represents the desired transmission rate. 
Consequently, no resources are utilized if the UE experiences poor fading situation.

We make the assumption that each BS serves all of its $U$ associated UEs, 
with the necessary radio resources allocated to fulfill their requirements. 
Consequently, the resource consumption within each cell is determined 
by the total resources allocated to the UEs within that cell. 
Furthermore, we assume that all BSs remain active during network transmissions, 
thereby the average resource consumption is determined by the average consumption of a typical Voronoi cell $\mathcal{V}_0$:
%\vspace{-8 pt}
\begingroup\makeatletter\def\f@size{10.}\check@mathfonts
\begin{align}\label{EQ_WUC0}
W(\cdot) = \mathbb{E} \left\{ \sum_{k\in \mathcal{V}_0} w_k^{\rm UC} \right\}.
\end{align}	
\endgroup
Note that the expectation is with respect to $\Phi_u$ and $\Phi_b$, 
as well as with respect to the channel coefficients, based on \eqref{EQ_SINR}.
For the network-wide resource consumption, we then have:
\begin{thm}\label{Theorem_Wt}
		Assume an interference-limited frequency reuse 1 cellular network,
		wherein every UE is served by its nearest BS.
		The position of BSs and UEs adhere to PPPs with intensities $\lambda_b$ and
		$\lambda_u$, respectively. 
		BSs allocate bandwidth to UEs using SINR service threshold
		$\gamma$ as outlined in (\ref{EQ_serviceTh}).
		The average of required resources across whole network is then:
	\begingroup\makeatletter\def\f@size{9.}\check@mathfonts
	\begin{align}\label{EQ_WUC}
		W(\gamma) &= \frac{\lambda_u}{\lambda_b } \int_0^{w_{\rm th}}
		w \frac{d}{dw} \Big( \frac{\eta(w) - \gamma}{\eta(w) T(\eta(w)) - \gamma T(\gamma)} \Big) dw,
	\end{align}
	\endgroup
	where $w_{\rm th}=\dfrac{R}{\log_2(1+\gamma)}$
	and $T(x) = 1+\sqrt{x}\: \tan^{-1}(\sqrt{x})$ and $\eta(w) = 2^{\frac{R}{w}-1}$.
\end{thm}

\begin{proof}
First, we need to calculate the expected resources required by a typical User Equipment (UE) requesting services from the network.
Based on \eqref{EQ_SINR} and \eqref{EQ_serviceTh} and by defining $w_k = \dfrac{R}{\log(1+\gamma_k)}$, we get:
\begingroup\makeatletter\def\f@size{8.75}\check@mathfonts
\begin{align}\label{WQ:aux}
	&\mathbb{E}_{ {\gamma}_k }\left\{  \frac{R}{\log(1+\gamma_k)} \mathbbm{1}(\gamma_k \geq \gamma)\right\} =
	 \mathbb{E}_{ w_k }\bigg\{  \frac{R}{\log(1+\gamma_k)}  \mathbbm{1}(w_k \leq w_{\rm th})\bigg\} \notag \\
	&\qquad= \mathbb{E}_{w_k}\left\{ \frac{R}{\log(1+\gamma_k)} \:\big|\: w_k \leq w_{\rm th} \right\} \mbox{Pr}(w_k \leq w_{\rm th}) \notag \\
	&\qquad= \int_0^{w_{\rm th}} w \frac{d}{dw} \mbox{F}_{w_k}(w) dw \: \mbox{F}_{w_k}(w_{\rm th}) ,
\end{align}
\endgroup
where $w_{\rm th} = R/\log(1+\gamma)$, 
and $\mbox{F}_{w_k}(w)$ is the cumulative density function (CDF) of r.v. $w_i$ evaluated at $w$,
Based on \eqref{EQ_SINR}, we then obtain:
\begin{align*}
	\mbox{F}_{w_k}(w) =&~ 1 - \mathbb{P}\left\{ |h_0^k|^2\|\boldsymbol{x}_k\|^{-e}\leq I_k \eta(w) \right\} \\
	\overset{(a)}=&~  \mathbb{E}_{I_k}\big\{ \exp(- \eta(w) I_k \|\boldsymbol{x}_k \|^{e}) \big\},
\end{align*}
where $\eta(w) = 2^{\frac{R}{w}-1}$, and $(a)$ is obtained based on  $|h_0^k|^2 \sim \exp(1)$.
Considering that $I_k = \sum_{j \in \Phi_b\backslash\{0\}} |h_j^k|^2 \|\boldsymbol{x}_k-\boldsymbol{r}_k \|^{-e}$ 
and $|h_j^k|^2\sim\exp(1)$, we can get:
\begingroup\makeatletter\def\f@size{9}\check@mathfonts
\begin{align*}
	\mathbb{E}_{I_k}\big\{ \exp\left(- \eta(w) I_k \|\boldsymbol{x}_k \|^{e}\right) \big\} &=\!\!\! \prod_{j\in \Phi_b}\! \left( 1+\eta(w) \frac{\|\boldsymbol{x}_k-\boldsymbol{r}_j\|^{-e}}{\|\boldsymbol{x}_k\|^{-e}} \right)^{-1} \\
	&=\!\!\! \prod_{j\in \Phi_b}\! P_{w}\left(\frac{\|\boldsymbol{x}_k-\boldsymbol{r}_j \|^{-e}}{\|\boldsymbol{x}_k\|^{-e} }\right),
\end{align*}
\endgroup
where $~\: P_{w}(\zeta ) := \big( 1+\eta(w)\: \zeta  \big)^{-1}$.
Therefore, we get:
$$
\mbox{F}_{w_k}(w) = \prod_{j\in \Phi_b}\! P_{w}\left(\frac{\|\boldsymbol{x}_k-\boldsymbol{r}_j \|^{-e}}{\|\boldsymbol{x}_k\|^{-e} }\right).
$$
By using $\dfrac{d}{dw}\mbox{F}_{w_k}(w) = \mbox{F}_{w_k}(w) \dfrac{d}{dw}\:\log\left( \mbox{F}_{w_k}(w)\right)$, 
and defining 
\begingroup\makeatletter\def\f@size{9}\check@mathfonts
\begin{align*}
	S_{w}(\zeta) := - \frac{d\eta(w)}{dw} \dfrac{\zeta}{1+\eta(w)\:\zeta},
\end{align*}
\endgroup
we can obtain:
\begingroup\makeatletter\def\f@size{9}\check@mathfonts
\begin{align}\label{EQ_Exp_w}
	&\!\!\!\! \mathbb{E}_{ {\gamma}_k  }\left\{  \frac{R}{\log(1+\gamma_k)} \mathbbm{1}(\gamma_k \geq \gamma)\right\}  \notag\\
	 & ~=\int_0^{w_{\rm th}} w \sum_{j\in\Phi_b}\hspace{-5 pt} S_{w} \left(\frac{\|\boldsymbol{x}_k-\boldsymbol{r}_j \|^{-e}}{\|\boldsymbol{x}_k\|^{-e} }\right) 
	\hspace{-5 pt}\prod_{j\in \Phi_b}\hspace{-5 pt} P_{w}\left(\frac{\|\boldsymbol{x}_k-\boldsymbol{r}_j \|^{-e}}{\|\boldsymbol{x}_k\|^{-e} }\right) dw \notag\\
	&\qquad\qquad\times \prod_{j\in \Phi_b}\hspace{-5 pt} P_{w_{\rm th}}\left(\frac{\|\boldsymbol{x}_k-\boldsymbol{r}_j \|^{-e}}{\|\boldsymbol{x}_k\|^{-e} }\right) .
\end{align}
\endgroup
\begin{figure*}[t]
	\vspace{5 pt}
	\begingroup\makeatletter\def\f@size{9.}\check@mathfonts
	\begin{equation}
		\begin{aligned}\label{EQ_button}
			W(\gamma) =&   \int_0^{w_{\rm th}} w \:\mathbb{E}_{\Phi_u}\bigg\{ \sum_{k\in\Phi_u} \lambda_b\iint_{\mathbb{R}^2} Q_{w}\left(\frac{\|\boldsymbol{x}_k-\boldsymbol{r} \|^{-e}}{\|\boldsymbol{x}_k\|^{-e} }\right) \mathbbm{1}\big(\boldsymbol{x}_k\in M_0(\boldsymbol{r})\big)rdr\:d\theta \\ 
			&\qquad \times \exp\bigg( \lambda_b \iint_{\mathbb{R}^2} \Big( P_{w}\left(\frac{\|\boldsymbol{x}_k-\boldsymbol{r} \|^{-e}}{\|\boldsymbol{x}_k\|^{-e} }\right) 
			P_{w_{\rm th}}\left(\frac{\|\boldsymbol{x}_k-\boldsymbol{r} \|^{-e}}{\|\boldsymbol{x}_k\|^{-e} }\right)
			\mathbbm{1}\big(\boldsymbol{x}_k\in M_0(\boldsymbol{r})\big)-1\Big) rdr\:d\theta \bigg) \bigg\}dw,
		\end{aligned}
	\end{equation}
	\endgroup
	\hrule height 0.1pt depth 0pt width 7.15in \relax
\end{figure*}
We now present the main part of the proof of Theorem \ref{Theorem_Wt}.
For this, we compute $W(\gamma) = \mathbb{E}\left\{ \sum_{k \in \mathcal{V}_0} w_k \right\}$ from \eqref{EQ_WUC0}.
\textcolor{black}{According to Slivnyak-Make theorem} \cite{BK1}, 
we assume that the BS responsible of cell $\mathcal{V}_0$ is located at origin, 
i.e., $\boldsymbol{r}_0=\boldsymbol{0}$. 
Based on \eqref{EQ_SINR}, we have:
\begingroup\makeatletter\def\f@size{7.75}\check@mathfonts
\begin{equation}\label{EQ_wt1}
	\begin{aligned}
		\!\!\! &W(\gamma)\hspace{-3 pt} ~= \mathbb{E}_{\Phi_u,\Phi_b} \mathbb{E}_{ \boldsymbol{\gamma} } \left\{  \sum_{k \in \mathcal{V}_0} \frac{R}{\log(1+\gamma_k)} \mathbbm{1}(\gamma_k \geq \gamma) \right\} \\
		&~~= \mathbb{E}_{\Phi_u,\Phi_b} \left\{  \bsum_{k \in \Phi_u}  \mathbb{E}_{ \gamma_k }\left\{  \frac{R}{\log(1+\gamma_k)} \mathbbm{1}(\gamma_k \geq \gamma)\right\} \mathbbm{1}(\boldsymbol{x}_k \in \mathcal{V}_0)   \right\} \\
		&~~\overset{(a)}=  \mathbb{E}_{\Phi_u}\bigg\{\hspace{-1 pt} \bsum_{k \in \Phi_u}  \mathbb{E}_{\Phi_b}\Big\{\hspace{-2 pt} \int_0^{w_{\rm th}}\hspace{-6 pt} w\hspace{-7 pt}\bsum_{j\in \Phi_b \backslash \{0\}}\hspace{-8 pt} S_{w}\hspace{-3 pt}\left(\hspace{-2 pt}\frac{\|\boldsymbol{x}_k-\boldsymbol{r}_j \|^{-e}}{\|\boldsymbol{x}_k\|^{-e} }\hspace{-2 pt}\right)\hspace{-3 pt} \mathbbm{1}\big(\boldsymbol{x}_k \in M_0(\boldsymbol{r}_j)\big) 
		\hspace{-8 pt}\\
		&~~~~~ \times\hspace{-10 pt} 
		\prod_{j\in \Phi_b \backslash \{0\}}\hspace{-13 pt} P_{w}\hspace{-2 pt}\left(\hspace{-2 pt}\frac{\|\boldsymbol{x}_k-\boldsymbol{r}_j \|^{-e}}{\|\boldsymbol{x}_k\|^{-e} }\hspace{-2 pt}\right)\!\!
		P_{w_{\rm th}}\hspace{-3 pt}\left(\hspace{-2 pt}\frac{\|\boldsymbol{x}_k-\boldsymbol{r}_j \|^{-e}}{\|\boldsymbol{x}_k\|^{-e} }\hspace{-2 pt}\right)
		\hspace{-2 pt}\mathbbm{1}\big(\boldsymbol{x}_k \!\in\! M_0(\boldsymbol{r}_j)\big) dw \Big\} \hspace{-2 pt}\bigg\}
	\end{aligned}
\end{equation}
\endgroup
where for $(a)$ we use \eqref{EQ_Exp_w} \textcolor{black}{and Lemma \ref{Lemma2}}.
In \eqref{EQ_wt1} we are dealing with a reduced-palm process as $j\in\Phi_b \backslash \{0\}$.
Considering that the distribution of reduced-palm process is equal 
to the distribution of original process \cite{BK1}, 
and based on Lemma \ref{Lemma1}, we obtain \eqref{EQ_button} 
which is written in the top of next page,
where $Q_{w}(\cdot) := S_{w}(\cdot)P_{w}(\cdot)P_{w_{\rm th}}(\cdot)$.
By following Lemma~\ref{Lemma2} that
$$
	\big\{ \boldsymbol{r} \in \mathbb{R}^2 \: |\: \textbf{x} \in M_0(\boldsymbol{r}) \big\} = B_c(\boldsymbol{r};\boldsymbol{x},\|\boldsymbol{x}\|),
$$
we then get:
\begingroup\makeatletter\def\f@size{9}\check@mathfonts
\begin{align}\label{EQ_app1}
		W(\gamma) &=  \int_0^{w_{\rm th}}\! w \lambda_b \lambda_u\! \iint_{\mathbb{R}^2} \Bigg\{\! \iint_{B_c(\boldsymbol{r}; \boldsymbol{x},\|\boldsymbol{x}\| )} \hspace{-15pt} Q_{w}\left(\frac{\|\boldsymbol{x}-\boldsymbol{r} \|^{-e}}{\|\boldsymbol{x}\|^{-e} }\right) rdr\:d\theta \notag \\ 
		&\!\!\!\!\!\!\!\!\times\!\exp\Big(\lambda_b\!\! \iint_{B_c(\boldsymbol{r}; \boldsymbol{x},\|\boldsymbol{x}\| )} \hspace{-15pt} P_{w}\left(\frac{\|\boldsymbol{x}-\boldsymbol{r} \|^{-e}}{\|\boldsymbol{x}\|^{-e} }\right)\!\!
		P_{w_{\rm th}}\left(\frac{\|\boldsymbol{x}-\boldsymbol{r} \|^{-e}}{\|\boldsymbol{x}\|^{-e} }\right)
		\! rdrd\theta \notag \\
		&~~~~~-\! \lambda_b \! \iint_{\mathbb{R}^2} \! rdrd\theta  \Big)\! \Bigg\}xdxd\theta'dw,
\end{align}
\endgroup
where $x=\| \boldsymbol{x}\|$.
Note that \eqref{EQ_app1} obtained by considering the probability generating functional (PGFL) property of PPP $\Phi_u$ \cite{BK1}.
Now, we consider change of variables $\boldsymbol{\rho} = \boldsymbol{x}-\boldsymbol{r}$, 
and then $\boldsymbol{q} = z\boldsymbol{\rho}$ to obtain:
\begingroup\makeatletter\def\f@size{8.5}\check@mathfonts
	\begin{align*}
		W(\gamma) =&~  \hspace{-3pt} \int_0^{w_{\rm th}} \hspace{-5pt} w \lambda_b \lambda_u \iint_{\mathbb{R}^2} \hspace{-3pt} \Big\{ 2\pi \hspace{-3pt} \int_1^\infty \hspace{-7pt} Q_{w}\left( z^{-e}  \right) q^2 zdz \\
		&\times\!\exp\Big( 2\pi\lambda_b q^2 \big(\int_1^\infty \hspace{-7pt} 
		P_{w}\!\!\left( z^{-e} \right)\!\!
		P_{w_{\rm th}}\!\!\left( z^{-e} \right)\!
		zdz - \int_0^\infty \hspace{-5pt} zdz\big)  \Big) \Big\}qdqdw \\[4 pt]
		\overset{(a)}=&  \int_0^{w_{\rm th}}  \Big( \frac{\lambda_u}{2\lambda_b}w\int_1^\infty Q_{w}(z^{-e})zdz \:   \Omega_l^2(w) \Big)dw \\[5 pt]
		=&~ \frac{\lambda_u}{2\lambda_b} \int_0^{w_{\rm th}} w \frac{d}{dw}\:\Omega_w(w) dw,	
	\end{align*}
\endgroup
where 
\begingroup\makeatletter\def\f@size{9}\check@mathfonts
$$
\Omega_w(w) = \frac{1}{\frac12-\int_1^\infty \left(P_{w}(z^{-e})P_{w_{\rm th}}(z^{-e})-1\right)zdz },
$$ 
\endgroup
$\rho=\| \boldsymbol{\rho}\|$, $q=\| \boldsymbol{q}\|$, $\dfrac{d}{dw}P_{w}(\cdot)P_{w_{\rm th}}(\cdot)=Q_{w}(\cdot)$ 
and for $(a)$ we use $\int_0^\infty q^3 \exp(\phi_0 q^2) dq = \dfrac{1}{2\phi_0^2}$ for $\phi_0<0$.
By computing the inner integral for $e=4$, the statement follows.
\end{proof}

\noindent \textbf{Remark}:
	For the case of service threshold $\gamma \to 0$,
	the average of total resource consumption tends to infinity:
	\begingroup\makeatletter\def\f@size{9.}\check@mathfonts
	\begin{align*}
		W(0) \to \infty.
	\end{align*}
	\endgroup
	This shows the importance of considering a truncated policy as \eqref{EQ_serviceTh}.

\noindent \textbf{Comparison to an approximation}:
It is possible to obtain an approximation for the average amount of resource consumption, 
denoted by $W^{\rm app}$.
Given that all BSs remain active during network transmissions, 
the average resource consumption is obtained by considering a typical Voronoi cell $\mathcal{V}_0$.
As such, we need to compute the average of needed resources in the cell $\mathcal{V}_0$.
One conventional way is thus to compute the average bandwidth consumption for a single UE, 
and then multiply that with the average number of UEs located at $\mathcal{V}_0$.
According to \eqref{WQ:aux}, and for path-loss exponent $e=4$, 
the average resource consumption for a typical UE located at origin ($\textbf{x}_0=\textbf{0}$) is
\begingroup\makeatletter\def\f@size{8.75}\check@mathfonts
\begin{align*}
	\bar{w}_0= \int_0^{w_{\rm th}} w \frac{d}{dw} \mbox{F}_{w_0}(w) dw \: \mbox{F}_{w_0}(w_{\rm th}) ,
\end{align*}
\endgroup
where 
\begingroup\makeatletter\def\f@size{9}\check@mathfonts
\begin{align*}
	\mbox{F}_{w_0}(w) =&~ 1 - \mathbb{P}\left\{ |h_0^0|^2\|\boldsymbol{r}_0\|^{-4}\leq I_0 \eta(w) \right\} \\
	\overset{(a)}=&~  \mathbb{E}_{I_0,\:r_0}\big\{ \exp(- \eta(w) I_0 \|\boldsymbol{r}_0 \|^{4}) \big\}\\
	=&~ \mathbb{E}_{\Phi_b,\{h_j^0\}_j,r_0}\Big\{ \prod_{j\in\Phi_b\backslash\{0\}} \exp(-\eta(w) |h_j^0|^2 \:r_j^{-4} \: r_0^4 ) \Big\}\\
	\overset{(b)}=&~ \mathbb{E}_{\Phi_b,r_0}\Big\{ \prod_{j\in\Phi_b\backslash\{0\}} (1+\eta(w)  \:r_j^{-4} \:r_0^4 )^{-1} \Big\}\\
	\overset{(c)}=&~ \mathbb{E}_{\:r_0}\Big\{ \exp\big(2\pi\lambda_b \int_{r_0}^\infty \big((1+\eta(w)  \:r^{-4} \:r_0^4)^{-1}-1 \big) rdr\big) \Big\}\\
	=&~ \mathbb{E}_{\:r_0}\Big\{ \exp\big(-\pi\lambda_b\:r_0^2 \sqrt{\eta(w)}\: \tan^{-1}(\sqrt{\eta(w)})   \big) rdr\big) \Big\}\\
	\overset{(e)}=&~ \frac{1}{1 + \sqrt{\eta(w)}\tan^{-1}(\sqrt{\eta(w)})},
\end{align*}
\endgroup
with $I_0 = \sum_{j\in\Phi_b \backslash \{0\}} |h_j^0|^2 \|\boldsymbol{r}_j\|^{-4}$ 
being the interference power experienced by UE $0$.
Moreover, for (a) and (b), we used $|h_0^0|^2 \sim \exp(1)$ and $|h_j^0|^2 \sim \exp(1)$,
and we set $r_0 = \| \boldsymbol{r}_0\|$ and $r_j = \| \boldsymbol{r}_j\|$,
for (c), we leveraged Campbell’s Theorem \cite{Kingman93},
for (d), we exploited the PDF of $r_0$: $f(r_0)=2\pi\lambda_b r_0 \exp(-\pi\lambda_br_0^2)$.

On the other hand, based on \eqref{EQ:avg_n}, the average number of UEs in the cell $\mathcal{V}_0$ 
is $\bar{n}_u=\frac{\lambda_u}{\lambda_b}$.
Therefore, we can get:
\begingroup\makeatletter\def\f@size{9.}\check@mathfonts
\begin{align}\label{EQ_W_app}
	W^{\rm app}(\gamma) &= \bar{n}_u \: \bar{w}_0 \notag \\
	&= \frac{\lambda_u}{\lambda_b } \int_0^{w_{\rm th}}
	w \frac{d}{dw} \Big( \frac{1}{1 + \sqrt{\eta(w)}\tan^{-1}(\sqrt{\eta(w)}) } \Big) dw \notag \\
	&\qquad \times \frac{1}{1 + \sqrt{\gamma}\tan^{-1}(\sqrt{\gamma}) },
\end{align}
\endgroup
Therefore, equation \eqref{EQ_W_app} presents an approximation of exact value of network-wide average resource consumption \eqref{EQ_WUC}.
\noindent \textbf{Remark}:
While equations \eqref{EQ_WUC} and \eqref{EQ_W_app} demonstrate resource consumption proportionate to $\lambda_u/\lambda_b$, 
the dependency on $\gamma$ differs between them.

\begin{figure}[t]
	\centering
	\includegraphics[width=9 cm]{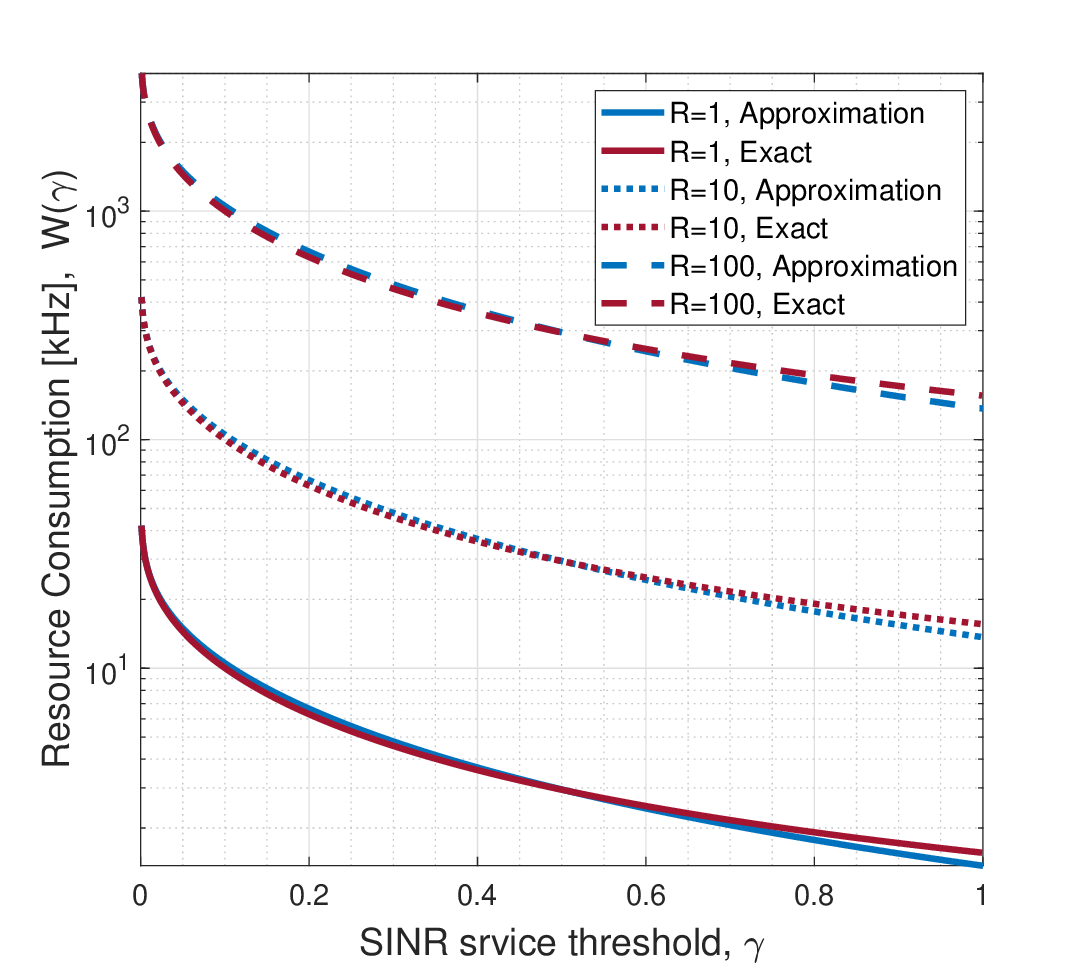}
	\caption{ Comparison between approximation $W^{\rm app}(\gamma)$ and \\exact $W(\gamma)$ values
		as a function of SINR service threshold $\gamma$. \label{Fig:App_exact}}
\end{figure}
\begin{figure}[t]		
	\includegraphics[width=9 cm]{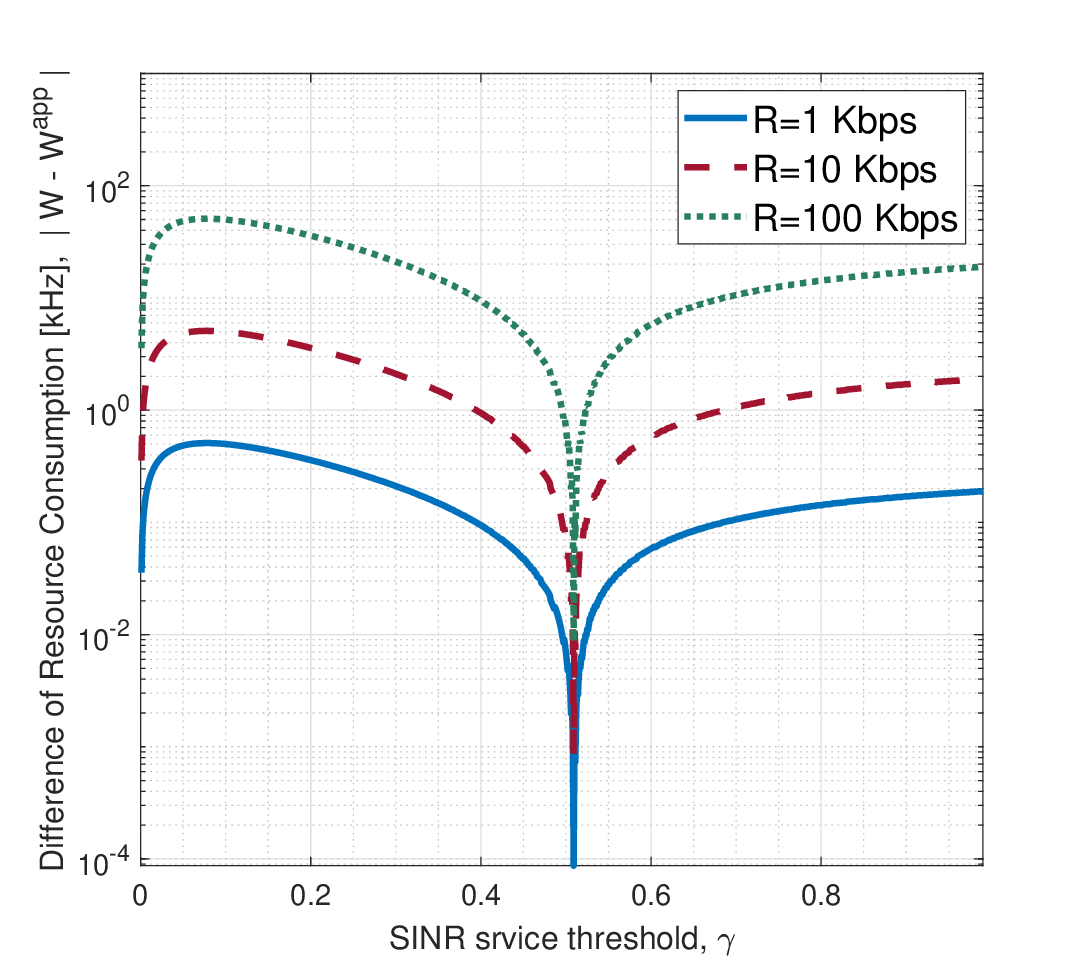}
	\caption{ Difference between approximation $W^{\rm app}(\gamma)$ and exact $W(\gamma)$ values
		as a function of SINR service threshold $\gamma$. \label{Fig:Diff}}
\end{figure}

Figures \ref{Fig:App_exact} and \ref{Fig:Diff} compare the approximation $W^{\rm app}$~\eqref{EQ_W_app} 
and the exact value $W$~\eqref{EQ_WUC} 
of network-wide resource consumption as a function of SINR service threshold $\gamma$ 
for $R \in \{1,10,100\}\: Kbps$ and ${\lambda_u}/{\lambda_b}=10$.
This result highlights the difference between $W^{\text{app}}$ and $W$, 
particularly for high and low values of $\gamma$. 
However, despite this variance, the approximation remains 
acceptable for moderate values of $\gamma$.

%/\/\/\/\/\/\/\/\/\/\/\/\/\/\/\/\/\/\/\/\/\/\/\/\/\/\/\/\/\/\/\/\/\/\/\/\/\/\/\/\/\/\/\/\/\/\/\/\/\/\/\/\/\/\/\/%
%  													   	 SubSec_Delivery										%
%/\/\/\/\/\/\/\/\/\/\/\/\/\/\/\/\/\/\/\/\/\/\/\/\/\/\/\/\/\/\/\/\/\/\/\/\/\/\/\/\/\/\/\/\/\/\/\/\/\/\/\/\/\/\/\/%

%/\/\/\/\/\/\/\/\/\/\/\/\/\/\/\/\/\/\/\/\/\/\/\/\/\/\/\/\/\/\/\/\/\/\/\/\/\/\/\/\/\/\/\/\/\/\/\/\/\/\/\/\/\/\/\/%
%  													   	 Conclusoin 											%
%/\/\/\/\/\/\/\/\/\/\/\/\/\/\/\/\/\/\/\/\/\/\/\/\/\/\/\/\/\/\/\/\/\/\/\/\/\/\/\/\/\/\/\/\/\/\/\/\/\/\/\/\/\/\/\/%
%\renewcommand{\baselinestretch}{1.08} 
\section{Conclusion}\label{Sec_Conclu}
In this paper, we developed two identities related to Poisson Point Processes and Poisson-Voronoi tessellation.
We then presented some applications of these identities.
One application was related to the computation of resource consumption 
for the conventional transmission scheme in cellular networks.
The results showed the versatility of derived identities in analysis of cellular networks.
Analysis of network performances for beamforming-based transmission schemes is considered as a future work.
%\section*{Acknowledgment}
%This work was funded by  Academy of Finland (grant 319058).
%The work of G. Caire was partially funded by  European Research Council under  ERC Advanced Grant N. 789190, CARENET.

%\small{
	\bibliographystyle{IEEEtran}
	\bibliography{IEEEabrv,IEEE}

\appendix
\subsection{Proof of Lemma~\ref{Lemma1}}\label{App1}
Let $N$ denote the number of points of $\Phi_1$, which follows a Poisson distribution.
Suppose that the points are placed in the region $A$ as a subspace of Cartesian space.
Note that for a homogeneous PPP, the points are independently and uniformly distributed 
within $A$ \cite{Chiu2013}. 
We then have:
\begingroup\makeatletter\def\f@size{7.75}\check@mathfonts
\begin{align*}
	&\mathbb{E}\Big\{ \prod_{k\in \Phi} P(\boldsymbol{x}_k) \sum_{k\in \Phi} S(\boldsymbol{x}_k) \Big \} =
	\mathbb{E}_N \bigg\{ \mathbb{E} \Big\{\sum_{k=1}^N S(\boldsymbol{x}_k)P(\boldsymbol{x}_k)\prod_{i\ne k}^N P(\boldsymbol{x}_i) \Big| N\Big\}\!\bigg\}\\
	&~\overset{(a)}= \mathbb{E}_N\left\{ \sum_{n=1}^N \iint_A \frac{S(\boldsymbol{s})P(\boldsymbol{s})}{|A|}d\boldsymbol{s} \: \Big(\iint_A \frac{P(\boldsymbol{s})}{|A|} d\boldsymbol{s} \Big)^{N-1}  \right\} \\
	&~= \iint_A \frac{S(\boldsymbol{s})P(\boldsymbol{s})}{|A|}d\boldsymbol{s} ~\mathbb{E}_N\left\{ N  \: \Big(\iint_A \frac{P(\boldsymbol{s})}{|A|} d\boldsymbol{s}\Big)^{N-1} \right\} \\
	&~\overset{(b)}=  \iint_A \frac{S(\boldsymbol{s})P(\boldsymbol{s})}{|A|}d\boldsymbol{s} \: \lambda_1 |A| \: \exp\big(-\lambda_1|A| \big) \: \exp\left(\lambda_1 |A| \iint_A \frac{P(\boldsymbol{s})}{|A|}d\boldsymbol{s} \right)
\end{align*}
\endgroup
where $|A|$ is the area of $A$, for $(a)$ we consider that points are independently and uniformly distributed over $A$,
and for $(b)$ we use $N\sim \mbox{Pois}(\lambda_1 |A|)$.
By rearranging the last equation and letting
$A=\mathbb{R}^2$, the statement follows.

\subsection{Proof of Lemma~\ref{Lemma2}}\label{App2}
\begin{figure}[t]
 	\vspace{10 pt}
	\centering
	\includegraphics[width=8.5 cm]{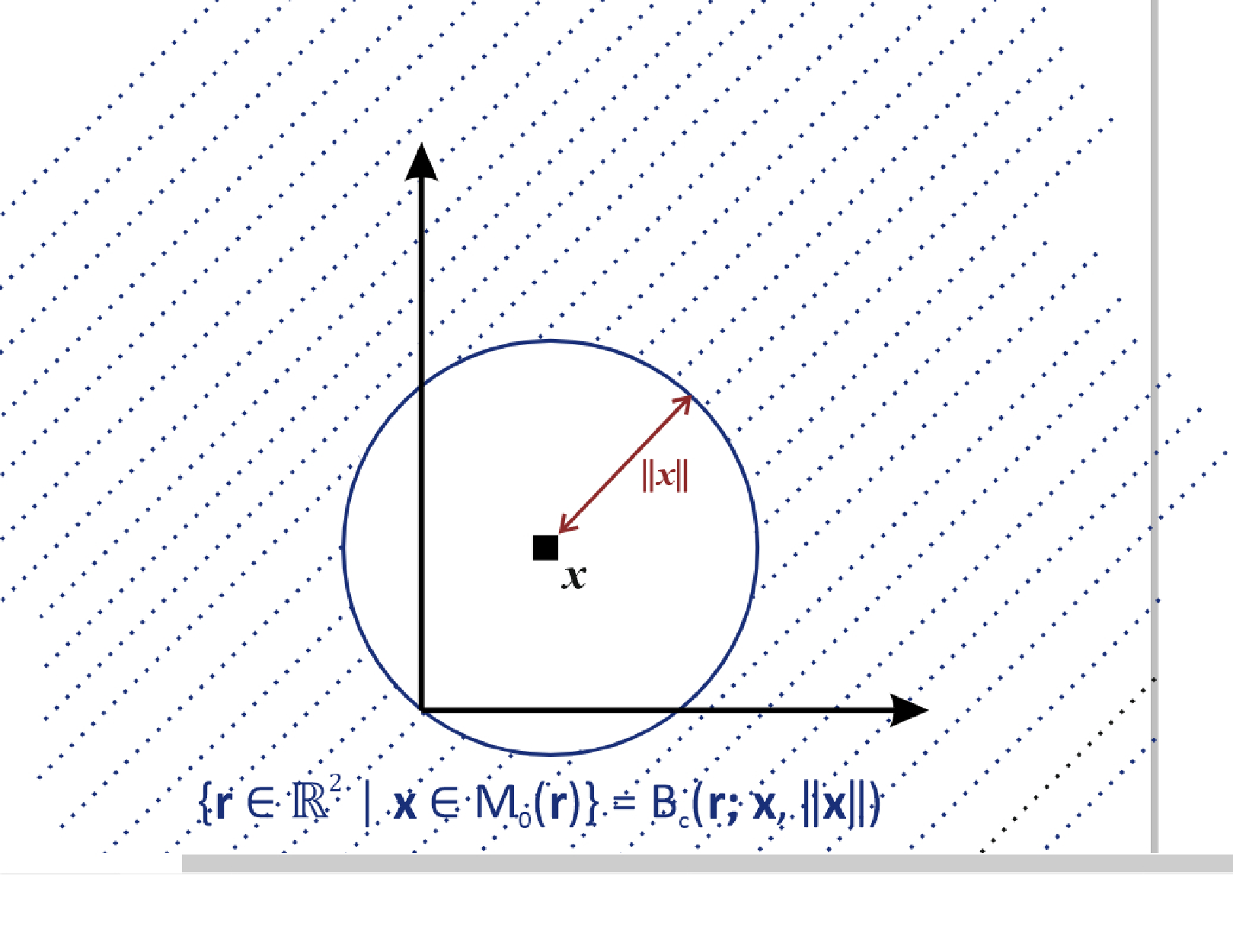}\vspace{10pt}
	\caption{ Illustration of $B_c(\textbf{r};\textbf{x},\|\textbf{x}\|)$. \label{Fig:Bc}}
\end{figure}

	\textcolor{black}{
	The first statement~\ref{EQ:State1} follows based on the geometrical definition of the Voronoi tessellation of $\Phi_2$
	and definition of $M_0(\cdot)$. 
	For the second statement~\ref{EQ:State2}, we consider the given $\textbf{x} \in \mathbb{R}^2$.
	Then, we exclude from the Euclidean space the area where $\textbf{x} \notin M_0(\textbf{r})$.
	The resulting area will be equal to a circle centered at $\textbf{x}$ with radius $\|\textbf{x}\|$,
	which is the definition of $B_c(\boldsymbol{r};\textbf{x},\|\textbf{x}\|)$.
	See Figure \ref{Fig:Bc} for an illustration.}

\subsection{Proof of Remark~\ref{Remark}}\label{App3}
\begingroup\makeatletter\def\f@size{8.75}\check@mathfonts
\begin{align*} 
	&\mathbb{E}\left\{ \sum_{k\in \mathcal{V}_0} P(\textbf{x}_k) \right\} 
	= \mathbb{E}\left\{ \sum_{k\in\Phi_1 }  P(\textbf{x}_k)  \mathbbm{1}(\textbf{x}_k \in \mathcal{V}_0 ) \right\} \\
	&\quad= \mathbb{E}_{\Phi_1,\:\Phi_2}\left\{ \sum_{k\in\Phi_1 } \prod_{j \in \Phi_2\backslash\{0\}}  P(\textbf{x}_k)  \mathbbm{1}(\textbf{x}_k \in M_0(\textbf{r}_j) ) \right\} \\
	&\quad= \mathbb{E}\left\{ \sum_{k\in\Phi_1 } P(\textbf{x}_k) \exp\Big(\lambda_2 \!\iint \big(\mathbbm{1}(\textbf{x}_k \in M_0(\textbf{r})-1 \big) d\textbf{r}  \Big) \right\} \\
	&\quad= \mathbb{E}\left\{ \sum_{k\in\Phi_1 } P(\textbf{x}_k) \exp\Big(\lambda_2 \!\iint_{B_c(\textbf{r},\textbf{x}_k,\|\textbf{x}_k\|)} \hspace{-25 pt} 1\: d\textbf{r} - \lambda_2 \!\iint_{\mathbb{R}^2}\!\! 1\: d\textbf{r} \Big) \right\} \\
	&\quad= \mathbb{E}\left\{ \sum_{k\in\Phi_1 } P(\textbf{x}_k) \exp\Big(-\pi\lambda_2  \|\textbf{x}_k\|^2 \Big) \right\} \\
	&\quad= \lambda_1 \iint P(\textbf{x})\exp\Big(-\pi\lambda_2  \|\textbf{x}\|^2 \Big) d\textbf{x},
\end{align*}
\endgroup
where obtained based on the Campbell theorem,
and considering that the distribution of reduced-palm process (when dealing with $\Phi_2\backslash\{0\}$) is equal  to the distribution of original process.
In a similar way, it can be shown that:
\begingroup\makeatletter\def\f@size{8.5}\check@mathfonts
\begin{align*} 
	&\mathbb{E}\left\{ \prod_{k\in \mathcal{V}_0} P(\textbf{x}_k) \right\} 
	= \exp\left( \lambda_1 \!\!\iint_{\mathbb{R}^2} \!\!\Big(P(\textbf{x})\exp\Big(-\pi\lambda_2  \|\textbf{x}\|^2 \Big) -1 \Big) d\textbf{x} \right).
\end{align*}
\endgroup

\end{document}